\documentstyle[12pt]{article}
\newcommand{\sect}[1]{\section{#1}\setcounter{equation}{0}}

\font\mbn=msbm10 scaled \magstep1
\font\mbs=msbm7 scaled \magstep1
\font\mbss=msbm5 scaled \magstep1
\newfam\mbff
\textfont\mbff=\mbn
\scriptfont\mbff=\mbs
\scriptscriptfont\mbff=\mbss\def\mbf{\fam\mbff}
\def\Re{{\mbf R}}

\def\Z{{\mbf Z}}
\def\Co{{\mbf C}}

\def\N{{\mbf N}}

\newtheorem{Th}{Theorem}[section]
\newtheorem{Lm}[Th]{Lemma}
\newtheorem{C}[Th]{Corollary}
\newtheorem{D}[Th]{Definition}
\newtheorem{Proposition}[Th]{Proposition}
\newtheorem{R}[Th]{Remark}
\newtheorem{E}[Th]{Example}
\author{Alexander Brudnyi\thanks{Research supported in part by NSERC.
}\\
Department of Mathematics and Statistics\\
University of Calgary, Calgary\\
Canada}
\title{Holomorphic Functions of Slow Growth on Coverings of
Pseudoconvex Domains in Stein Manifolds}
\date{} 
\begin{document} 
\maketitle
\begin{abstract}
{We apply the methods developed in [Br1] 
to study holomorphic functions of slow growth on coverings of pseudoconvex
domains in Stein manifolds. In particular,
we extend and strengthen certain results of 
Gromov, Henkin and Shubin [GHS] on holomorphic $L^{2}$ functions on
coverings of pseudoconvex manifolds in the case of coverings of Stein
manifolds.}
\end{abstract}

\sect{\hspace*{-1em}. Introduction.}
Let $M$ be a complex manifold satisfying 
\begin{equation}\label{e1}
M\subset\subset\widetilde M\subset N\ \ \ {\rm and}\ \ \
\pi_{1}(M)=\pi_{1}(N)
\end{equation}
where $M$ and $\widetilde M$ are open connected subsets of a complex
manifold $N$ and $\widetilde M$ is Stein. (Here $\pi_{1}(X)$ stands for
the fundamental group of $X$.) Condition (\ref{e1}) is valid, e.g., for
$M$ a strictly pseudoconvex domain or an analytic polyhedra in a 
Stein manifold. It implies that the group 
$\pi_{1}(N)$ is finitely generated. In [Br1] we presented a method to
construct integral representation
formulas for holomorphic functions of slow growth defined on unbranched
coverings of $M$. Using such formulas we established that some known
results for holomorphic functions on $M$ can be extended to similar results
for holomorphic functions of slow growth on coverings of $M$. In this paper
we continue to study holomorphic functions of slow growth
on coverings of $M$ and
apply the methods developed in [Br1] to extend and strengthen certain 
results of Gromov, Henkin and Shubin [GHS] on holomorphic $L^{2}$ functions 
on coverings of
pseudoconvex manifolds in the case of coverings of Stein manifolds.

The presentation in this paper is focused on several problems
and results formulated in [GHS]. To describe them we, first, recall
some definitions.

Let $M\subset\subset N$ be a domain with a smooth boundary 
$bM$ in an $n$-dimensional complex manifold $N$, that is,
\begin{equation}\label{m1}
M=\{z\in N\ :\ \rho(z)<0\}
\end{equation}
where $\rho$ is a real-valued function of class $C^{2}(\Omega)$ in a
neighbourhood $\Omega$ of the compact set $\overline{M}:=M\cup bM$
such that
\begin{equation}\label{m2}
d\rho(z)\neq 0\ \ \ {\rm for\ all}\ \ \ z\in bM\ .
\end{equation}
Let $z_{1},\dots, z_{n}$ be complex local coordinates in $N$ near $z\in bM$.
Then the tangent space $T_{z}N$ at $z$ is identified with $\Co^{n}$.
By $T_{z}^{c}(bM)\subset T_{z}N$ we denote the complex tangent space to
$bM$ at $z$, i.e.,
\begin{equation}\label{m3}
T_{z}^{c}(bM)=\{w=(w_{1},\dots,w_{n})\in T_{z}(N)\ :\ \sum_{j=1}^{n}
\frac{\partial\rho}{\partial z_{j}}(z)w_{j}=0\}\ .
\end{equation}
The {\em Levi form} of $\rho$ at $z\in bM$ is a hermitian form on
$T_{z}^{c}(bM)$ defined in the local coordinates by the formula
\begin{equation}\label{m4}
L_{z}(w,\overline{w})=\sum_{j,k=1}^{n}
\frac{\partial^{2}\rho}{\partial z_{j}\partial\overline{z}_{k}}(z)w_{j}
\overline{w}_{k}\ .
\end{equation}
The manifold $M$ is called {\em pseudoconvex} if $L_{z}(w,\overline{w})\geq 0$
for all $z\in bM$ and $w\in T_{z}^{c}(bM)$. It is called {\em strictly
pseudoconvex} if $L_{z}(w,\overline{w})>0$ for all $z\in bM$ and all
$w\neq 0$, $w\in T_{z}^{c}(bM)$.

Equivalently, strictly pseudoconvex manifolds can be described as the ones 
which locally, in a neighbourhood of any boundary point, can be presented as 
strictly convex domains in $\Co^{n}$. It is also known (see [C], [R]) that 
any strictly pseudoconvex manifold admits a proper holomorphic map with 
connected fibres onto a normal Stein space.

Without loss of generality we may and will assume that 
$\pi_{1}(M)=\pi_{1}(N)$ for $M$ as above.
Let $r: N_{G}\to N$ be the regular covering of $N$ with (discrete) 
transformation group $G$. Then $M_{G}:=r^{-1}(M)$ is a regular covering of
$M$ (with the same transformation group). It is a domain in $N_{G}$ with
the smooth boundary $bM_{G}:=r^{-1}(bM)$. By
$\overline{M}_{G}:=M_{G}\cup bM_{G}$ we denote the closure of $M_{G}$ in
$N_{G}$.

Let  $X$ be a subspace of the space 
${\cal O}(M_{G})$ of all holomorphic functions on $M_{G}$.
A point $z\in bM_{G}$ is called a {\em peak point} for $X$ if there exists
a function $f\in X$ such that $f$ is unbounded on $M_{G}$ but bounded outside
$U\cap M_{G}$ for any neighbourhood $U$ of $z$ in $N_{G}$.

A point $z\in bM_{G}$ is called a {\em local peak point} for $X$ if there
exists a function $f\in X$ such that $f$ is unbounded in $U\cap M_{G}$ for
any neighbourhood $U$ of $z$ in $N_{G}$ and there exists a 
neighbourhood $U$ of $z$ in $N_{G}$ such that for any neighbourhood $V$ 
of $z$ in $N_{G}$ the function $f$ is bounded on $U\setminus V$.

The Oka-Grauert theorem [G1] states that if $M$ is strictly pseudoconvex 
and $bM$ is not empty then every $z\in bM$ is a peak point for ${\cal O}(M)$.
In general it is not known whether the similar statement is true for 
boundary points of $M_{G}$ with an infinite $G$.

Let $dV_{M_{G}}$ be the Riemannian volume form on $M_{G}$ obtained by
a Riemannian metric pulled back from $N$. By $H^{2}(M_{G})$ we denote the
Hilbert space of holomorphic functions $g$ on $M_{G}$ with norm
$$
\left(\int_{z\in M_{G}}|g(z)|^{2}dV_{M_{G}}(z)\right)
^{1/2} .
$$

In [GHS] the von Neumann $G$-dimension $dim_{G}$ was used to measure the
space $H^{2}(M_{G})$. In particular, in [GHS, Theorem 0.2] the 
following result was proved.\\
{\bf Theorem A}\ {\em If $M$ is strictly pseudoconvex, then}
\begin{itemize}
\item[(a)]
$dim_{G}H^{2}(M_{G})=\infty$\ \ \ {\em and}
\item[(b)]
{\em each point in $bM_{G}$ is a local peak point for $H^{2}(M_{G})$.}
\end{itemize}

In [GHS, Theorem 0.5] the similar result was established
for a covering $M_{G}$ of a pseudoconvex manifold $M$ with a strictly 
plurisubharmonic $G$-invariant function existing in a neighbourhood of
$bM_{G}$. Finally, in [GHS, section 4] the following open problems were
formulated.

Suppose that $M$ is strictly pseudoconvex.

{\bf 1.} Does there exist a finite number of functions in 
$H^{2}(M_{G})\cap C(\overline{M}_{G})$ which separate all points in 
$bM_{G}$?

{\bf 2.} Assume that $dim_{\Co}M=2$. Does there exist 
$f\in H^{2}(M_{G})\cap C(\overline{M}_{G})$ such that $f(x)\neq 0$ for all
$x\in bM_{G}$?

{\bf 3.} Is it true that for every $CR$-function 
$f\in L^{2}(bM_{G})\cap C(bM_{G})$ in case $dim_{\Co}M_{G}>1$
there exists $f'\in H^{2}(M_{G})\cap C(\overline{M}_{G})$ such that
$f'|_{bM_{G}}=f$? 

Here $L^{2}(bM_{G})$ is defined similar to 
$H^{2}(M_{G})$ with respect to the volume form on $bM_{G}$ obtained by a
Riemannian metric pulled back from $N$. Also, recall that
$f\in C(bM_{G})$ is called a $CR$-function if for every smooth 
$(n,n-2)$-form $\omega$ with a compact support one has
$$
\int_{bM_{G}}f\wedge\overline{\partial}\omega=0\ .
$$
If $f$ is smooth this is equivalent to the fact that $f$ is a solution of
the tangential $CR$-equations:
$\overline\partial_{b}f=0$ (see, e.g., [KR]).

The present paper deals with the
above results and problems in the case of coverings of $M$ satisfying
condition (\ref{e1}).
\sect{\hspace*{-1em}. Formulation of Main Results.}
{\bf 2.1.} We start with some results related to Problem 2 of
the Introduction.

Let $M$ be a manifold satisfying condition (\ref{e1}) and $M'$ be
an unbranched covering of $M$. Condition
(\ref{e1}) implies that there is a covering $r:N'\to N$ of $N$
such that $M'$ is a domain in $N'$ (i.e., $\pi_{1}(M')=\pi_{1}(N')$).
As above $\overline{M'}$ denotes the closure of $M'$ in $N'$.

Let $\phi:N'\to\Re$ be a function
uniformly continuous with respect to the path metric induced by a
Riemannian metric pulled back from $N$. 
\begin{Th}\label{te1}
There exist a function $f_{\phi}\in {\cal O}(M')\cap C(\overline{M'})$
and a constant $C=C(\phi,M',N')$\footnote{Here and below the notation
$C=C(\alpha,\beta,\gamma,\dots)$ means that the constant depends only on
the parameters $\alpha,\beta,\gamma,\dots$ .} such that 
$$
|f_{\phi}(z)-\phi(z)|<C\ \ \ {\rm and}\ \ \
|df_{\phi}(z)|<C\ \ \ {\rm for\ all}\ \ \ z\in M'\ .
$$
\end{Th}
(Here the norm $|\omega(z)|$ of a differential form $\omega$ 
at $z\in M'$ is determined with respect to the Riemannian metric pulled back
from $N$.)

As a corollary of this result we answer an extended version of Question 2
of the Introduction for coverings of manifolds $M$ satisfying
condition (\ref{e1}). 
Namely, let $d$ be the path metric on $N'$
obtained by the pullback of a Riemannian metric defined on $N$. Fix a point
$o\in M'$ and set
$$
d_{o}(x):=d(o,x)\ ,\ \ \ x\in N'\ .
$$
From the triangle inequality it follows that the function
$\phi(x):=d_{o}(x)$, $x\in N'$, satisfies the hypothesis of Theorem
\ref{te1}.
\begin{C}\label{co1}
Let $f:=f_{\phi}$ be the function from Theorem \ref{te1} for $\phi=d_{o}$.
Then there exists a constant $\alpha>0$ such that 
$F=e^{-\alpha f}\in H^{2}(M')\cap C(\overline{M'})$.
\end{C}
(Note that here $F(x)\neq 0$ for all $x\in M'$ and
there are no restrictions on $dim_{\Co}M$.)\\
{\bf 2.2.} In this part we formulate our results related to Theorem A of the
Introduction. 

Let $M'$ be an unbranched covering of $M$ satisfying (\ref{e1}).
Let $\psi:M'\to\Re_{+}$ be a continuous function and $dV_{M'}$ be the 
Riemannian volume form on $M'$ obtained by a Riemannian metric 
pulled back from $N$. For an open set $D\subset M$ we introduce the Banach 
space $H^{p}_{\psi}(D')$ of holomorphic functions $g$ on 
$D':=r^{-1}(D)\subset M'$ with norm
$$
\left(\int_{z\in M'}|g(z)|^{p}\psi(z)dV_{M'}(z)\right)^{1/p} .
$$

Let $r:N'\to N$ be the covering of $N$ satisfying (\ref{e1}) such that
$\pi_{1}(N')=\pi_{1}(M')$. Then $M'\ (=r^{-1}(M))$ is a domain in $N'$.
Suppose that $\psi:N'\to\Re_{+}$ is such that $\log\psi$ 
is uniformly continuous with respect to the path metric 
induced by a Riemannian metric pulled back from $N$. 
We set $\phi:=\log\psi$ and consider the holomorphic function $f_{\phi}$ from 
Theorem \ref{te1}. This theorem implies that for $\widetilde C:=e^{C}$ 
\begin{equation}\label{equiv1}
\frac{1}{\widetilde C}\psi(z)\leq |e^{f_{\phi}(z)}|\leq\widetilde C\psi(z)\ ,
\ \ \ z\in\overline{M'}\ .
\end{equation}
Therefore the following result holds.
\begin{Proposition}\label{pro1}
For any open set $D\subset M$ and every $p\in [1,\infty)$ the map
$L_{\psi}: H_{\psi}^{p}(D')\to H_{1}^{p}(D')$, 
$L_{\psi}(g)=g\cdot e^{f_{\phi}/p}$, is an isomorphism of Banach 
spaces. \ \ \ \ \ \ \ \ $\Box$
\end{Proposition}

Let us now formulate an extension of Theorem A of the Introduction.

Suppose that $M$ is a strictly pseudoconvex domain in a complex manifold
$N$ such that $\pi_{1}(M)=\pi_{1}(N)$ and $N$ is a domain in a Stein 
manifold.  Let $r:N'\to N$ be an unbranched covering
of $N$ and $M'=r^{-1}(M)$ be the corresponding covering of $M$. Let
$bM'=r^{-1}(bM)$ be the boundary of $M'$ in $N'$. 
\begin{Th}\label{te2}
Each point in $bM'$ is a peak point for ${\cal O}(M')$ and for every
$H_{1}^{p}(M')$, $1\leq p<\infty$.
\end{Th}

From Theorem \ref{te2} and Proposition \ref{pro1} we get
(for $\psi$ as in Proposition \ref{pro1})
\begin{C}\label{co11}
Each $z\in bM'$ is a peak point for $H_{\psi}^{p}(M')$, $1\leq p<\infty$.
\ \ \ \ \ $\Box$
\end{C}
\begin{R}
{\rm The main ingredient of the proof of Theorem \ref{te2} is uniform
estimates for solutions of certain $\overline{\partial}$-equations on $M'$.
In fact similar estimates are valid on coverings of, so-called,
non-degenerate pseudoconvex polyhedrons on Stein manifolds 
(see [SH] and [He] for their definition). This class contains, in particular,
piecewise strictly pseudoconvex domains and non-degenerate
analytic polyhedrons on Stein manifolds. Also, every $M$ from this
class satisfies condition (\ref{e1}). Let
$M'$ be a covering of such $M$ and $z\in bM'$ be such that $M'\cap U$ is
strictly pseudoconvex for a neighbourhood $U\subset N'$ of $z$. Then arguing
as in the proof of Theorem \ref{te2} one obtains that {\em $z$ is a peak point
for ${\cal O}(M')$ and for every $H_{1}^{p}(M')$, $1\leq p<\infty$.}}
\end{R}
{\bf 2.3.} In this section we discuss some results related to Problem 3 of
the Introduction.

Let $r:N'\to N$ be a covering of $N$ satisfying (\ref{e1}). As before we
set $M'=r^{-1}(M)\subset N'$. Consider a function
$\psi: N'\to\Re_{+}$ such that
$\log\psi$ {\em is uniformly continuous with respect to the path metric 
induced by a Riemannian metric pulled back from} $N$. 
For such $\psi$ and every $x\in M$ we introduce the 
Banach space $l_{p,\psi, x}(M')$, $1\leq p\leq\infty$,
of functions $g$ on $r^{-1}(x)\subset M'$ with norm
\begin{equation}\label{e5}
|g|_{p,\psi,x}:=\left(\sum_{y\in
r^{-1}(x)}|g(y)|^{p}\psi(y)\right)^{1/p} .
\end{equation}
Next, for an open set $D\subset M$ we 
introduce the Banach space ${\cal H}_{p,\psi}(D')$, $1\leq p\leq\infty$, of 
functions $f$ holomorphic on $D':=r^{-1}(D)\subset M'$ with norm
\begin{equation}\label{e3}
|f|_{p,\psi}^{D}:=\sup_{x\in D}|f|_{p,\psi,x}\ .
\end{equation}

Clearly, one has a continuous embedding 
${\cal H}_{p,\psi}(D')\hookrightarrow H^{p}_{\psi}(D')$. 
Let $U\subset N$ be an open set
containing $\overline{D}$ and $U'=r^{-1}(U)$. Then for $\psi$ as 
above using the mean value property for plurisubharmonic functions
one can easily show that {\em
for each $p\in [1,\infty]$ the restriction $f\mapsto f|_{D'}$ induces a 
linear continuous map $H_{\psi}^{p}(U')\to {\cal H}_{p,\psi}(D')$.}
Also, for such $\psi$ from the results proved in [Br1] follow
that {\em holomorphic functions from ${\cal H}_{p,\psi}(M')$ separate 
all points in $M'$ (for each $p\in [1,\infty]$).}

Let us formulate the main result of this section.

Let $M\subset\subset\widetilde M\subset N$ be manifolds satisfying condition 
(\ref{e1}) with $dim_{\Co}M\geq 2$. Suppose that
$D\subset\subset M$ is an open subset whose boundary $bD$ is a connected
$C^{k}$ submanifold of $M$ $(1\leq k\leq\infty)$. 
For a covering $r:N'\to N$ we set $D'=r^{-1}(D)$ and $bD'=r^{-1}(bD)$.
\begin{Th}\label{te4}
For every $CR$-function $f\in C^{s}(bD')$, $0\leq s\leq k$, satisfying
$$
f|_{r^{-1}(x)}\in l_{p,\psi,x}(M')\ \ \ {\rm for\ all}\ \ \ x\in D \ \ \
{\rm and}\ \ \
\sup_{x\in bD}|f|_{p,\psi,x}<\infty
$$
there exists a 
function $f'\in {\cal H}_{p,\psi}(D')\cap C^{s}(\overline{D'})$ such 
that $f'|_{bD'}=f$.
\end{Th}
\begin{R}\label{re4}
{\rm (1) The converse to this theorem is always true: the restriction of
every $f'\in {\cal H}_{p,\psi}(D')\cap C^{s}(\overline{D'})$ to $bD'$ is
a $CR$-function satisfying the hypotheses of the theorem. \\
(2) We will prove also (see (\ref{norm3})) that for some $c=c(M',M,\psi,p)$}
$$
|f'|_{p,\psi}^{D}\leq c\sup_{x\in bD}|f|_{p,\psi,x}\ .
$$
\end{R}

As a corollary we obtain an analog of the Hartogs 
extension theorem. We formulate it for functions of the maximal possible 
growth for which our method works. 

Suppose that $M$ satisfies (\ref{e1}). Let
$D\subset M$ be a domain and $K\subset\subset D$ be a compact set such that
$U:=D\setminus K$ is connected. Consider a covering $r: M'\to M$
and set $D'=r^{-1}(D)$, $U'=r^{-1}(U)$. By $d_{o}$, $o\in M'$, we
denote the distance function on $M'$ as in Corollary \ref{co1}.
\begin{C}\label{co2}
There exists a constant $c>0$ such that for every 
$f\in {\cal O}(U')$ satisfying for some $c_{2}>0$ and $0<c_{1}<c$
the inequality
$$
|f(z)|\leq e^{c_{2}e^{c_{1}d_{o}(z)}}\ ,\ \ \ z\in U'\ ,
$$
there is $f'\in {\cal O}(D')$ such that
$$
|f'(z)|\leq e^{c_{3}e^{c_{1}d_{0}(z)}}\ ,\ \ \ z\in D'\ ,\ \ \ {\rm and}
\ \ \ f'|_{U'}=f\ ;
$$ 
where $c_{3}$ depends on $c_{2}, c_{1}, c, M, M'$ only.
\end{C}

We don't know whether a similar extension result holds for functions
$f$ growing faster than those of the corollary.\\
{\bf 2.4.} Finally, we formulate a result related to Problem 1 of the
Introduction.  
First, we recall some definitions of the theory of flat vector bundles
(see, e.g., [O]).

Let $X$ be a complex manifold and $\rho:\pi_{1}(X)\to GL_{k}(\Co)$ be
a homomorphism of its fundamental group. We set $G:=\pi_{1}(X)/Ker\ \!\rho$.
It is well known (see, e.g., Example \ref{constbun}
(b) below) that to any such $\rho$ corresponds a complex flat vector bundle
$E_{\rho}$ on $X$ (i.e., a bundle constructed by a locally constant
cocycle). We call $E_{\rho}$ the {\em bundle associated with} $\rho$. 
Assume that $\rho$ is such that $E_{\rho}$ is topologically trivial, i.e.,
is isomorphic in the category of continuous bundles to the bundle
$X\times\Co^{k}$. Every such $\rho$ can be obtained as the monodromy of the
equation $dF=\omega F$ on $X$ where $\omega$ is a matrix-valued $1$-form on
$X$ satisfying $d\omega-\omega\wedge\omega=0$.
By ${\cal T}(X)$ we denote the class of quotient groups $G$ obtained by
representations $\rho$ as above.

Now, let $r:M_{G}\to M$ be the regular covering of $M$ satisfying
condition (\ref{e1}) with transformation group $G$. Let $G_{1}\subset G$
be a subgroup of a finite index. Then there is a finite covering
$r_{1}:M_{1}\to M$ whose fibre is the quotient set $G/G_{1}$ such that
$M_{G}$ is also the regular covering of $M_{1}$ with transformation group
$G_{1}$.
\begin{Th}\label{te5}
Assume that $G_{1}\in {\cal T}(M_{1})$.
Then there is a finite number of functions
in ${\cal H}_{2,1}(M_{G})\cap C(\overline{M}_{G})$ which separate all points
in $M_{G}$.
\end{Th}
\begin{R}\label{re5}
{\rm (1) We will see from the proof that the functions in Theorem \ref{te5}
can be taken even from ${\cal H}_{2,\psi}(M_{G})$ where 
$\psi:M_{G}\to\Re_{+}$ has a double exponential growth.\\
(2) As the group $G$ in Theorem \ref{te5} one can take, e.g., 
a finitely generated free group (see, e.g., [O]) or a polycyclic group
(see, e.g., [Ra]). If $dim_{\Co}M=1$ then, since $M$ is homotopically 
equivalent to a 
one-dimensional $CW$-complex (see, e.g., [GR]), every quotient group $G$ 
obtained by a linear representation $\rho$ belongs to ${\cal T}(M)$.}
\end{R}
\sect{\hspace*{-1em}. Preliminary Results.}
{\bf 3.1.} First, we recall some basic facts from the theory of bundles
see, e.g., [Hi].

Let $X$ be a complex analytic space and $S$ be a complex 
analytic Lie group with
unit $e\in S$. Consider an effective holomorphic action of $S$ on
a complex analytic space $F$. Here {\em holomorphic action} means a 
holomorphic map $S\times F\rightarrow F$ sending 
$s\times f\in S\times F$ to
$sf\in F$ such that $s_{1}(s_{2}f)=(s_{1}s_{2})f$ and $ef=f$ for any
$f\in F$. {\em Efficiency} means that the condition $sf=f$ for some $s$ and 
any $f$ implies that $s=e$. 
\begin{D}\label{de2}
A complex analytic space $W$ together with a holomorphic map (projection)
$\pi:W\rightarrow X$ is called a holomorphic bundle on $X$ with 
structure group $S$ and fibre $F$, if there exists a
system of coordinate transformations, i.e., if
\begin{itemize}
\item[{\rm (1)}]
there is an open cover ${\cal U}=\{U_{i}\}_{i\in I}$ of $X$ and a
family of biholomorphisms 
$h_{i}:\pi^{-1}(U_{i})\rightarrow U_{i}\times F$,
that map ``fibres'' $\pi^{-1}(u)$ onto $u\times F$;
\item[{\rm (2)}] 
for any $i,j\in I$ there are elements 
$s_{ij}\in {\cal O}(U_{i}\cap U_{j}, S)$ such that
$$
(h_{i}h_{j}^{-1})(u\times f)=u\times s_{ij}(u)f\ \ \ {\rm for\ any}\ \
u\in U_{i}\cap U_{j},\ f\in F\ .
$$
\end{itemize}
A holomorphic bundle $\pi:W\rightarrow X$ whose fibre is a Banach space 
$F$ and the structure group is
$GL(F)$ (the group of linear invertible transformations of $F$) is
called a holomorphic Banach vector bundle. 
A holomorphic section of a holomorphic bundle $\pi:W\rightarrow X$
is a holomorphic map $s:X\rightarrow W$ satisfying $\pi\circ s=id$.
\end{D}
We will use the following construction of holomorphic bundles 
(see, e.g. [Hi,\ Ch.1]):

Let $S$ be a complex analytic Lie group and 
${\cal U}=\{U_{i}\}_{i\in I}$ be an open cover of $X$. 
By $Z_{\cal O}^{1}({\cal U},S)$ we denote the set of holomorphic $S$-valued 
${\cal U}$-cocycles. By definition, 
$s=\{s_{ij}\}\in Z_{\cal O}^{1}({\cal U},S)$, 
where $s_{ij}\in {\cal O}(U_{i}\cap U_{j}, S)$ and
$s_{ij}s_{jk}=s_{ik}\ \ \  {\rm on}\ \ \  U_{i}\cap U_{j}\cap U_{k}$.
Consider the disjoint union $\sqcup_{i\in I}U_{i}\times F$ and for any 
$u\in U_{i}\cap U_{j}$ identify the point $u\times f\in U_{j}\times F$ with 
$u\times s_{ij}(u)f\in U_{i}\times F$.
We obtain a holomorphic bundle $W_{s}$ on $X$ whose projection is induced 
by the projection $U_{i}\times F\rightarrow U_{i}$.
Moreover, any holomorphic bundle on $X$ with structure group $S$ and 
fibre $F$ is isomorphic (in the category of holomorphic bundles) to a 
bundle $W_{s}$. 
\begin{E}\label{constbun}
{\rm {\bf (a)}  Let $M$ be a complex manifold.
For any subgroup $H\subset\pi_{1}(M)$ consider the unbranched covering
$r:M(H)\rightarrow M$ corresponding to $H$. We will describe $M(H)$
as a holomorphic bundle on $M$.

First, assume that $H\subset\pi_{1}(M)$ is a normal subgroup.
Then $M(H)$ is a regular covering of $M$ and the quotient group
$G:=\pi_{1}(M)/H$ acts holomorphically on $M(H)$ by deck 
transformations.  It is well known that $M(H)$ in this case 
can be thought of as a {\em principle fibre bundle} on $M$ with fibre
$G$ (here $G$ is equipped with the discrete topology). 
Namely, let us consider the map $R_{G}(g):G\rightarrow G$, $g\in G$,
defined by the formula
$$
R_{G}(g)(q)=q\cdot g^{-1},\ \ \ q\in G.
$$ 
Then for an open cover ${\cal U}=\{U_{i}\}_{i\in I}$ of $M$ by sets
biholomorphic to open Euclidean balls in some $\Co^{n}$ there is
a locally constant cocycle $c=\{c_{ij}\}\in Z_{\cal O}^{1}({\cal U},G)$ 
such that $M(H)$ is biholomorphic
to the quotient space of the disjoint union 
$V=\sqcup_{i\in I}U_{i}\times G$ by the equivalence relation:
$U_{i}\times G\ni x\times R_{G}(c_{ij})(q)\sim x\times q\in 
U_{j}\times G$.
The identification space is a holomorphic bundle with projection 
$r:M(H)\rightarrow M$ induced by the
projections $U_{i}\times G\rightarrow U_{i}$.
In particular, when $H=e$ we obtain the definition of the universal
covering $M_{u}$ of $M$. 

Assume now that $H\subset\pi_{1}(M)$ is not necessarily normal.
Let $X_{H}=\pi_{1}(M)/H$ be the set of cosets with respect
to the (left) action of $H$ on $\pi_{1}(M)$ defined by left multiplications. 
By $[Hq]\in X_{H}$ we denote the coset containing $q\in\pi_{1}(M)$.
Let $A(X_{H})$ be the group of all homeomorphisms of $X_{H}$
(equipped with the discrete topology). We define the homomorphism
$\tau:\pi_{1}(M)\rightarrow A(X_{H})$ by the formula:
$$
\tau(g)([Hq]):=[Hqg^{-1}],\ \ \ q\in\pi_{1}(M).
$$
Set $Q(H):=\pi_{1}(M)/Ker(\tau)$ and let 
$\widetilde g$ be the image of $g\in\pi_{1}(M)$ in $Q(H)$. 
By $\tau_{Q(H)}:Q(H)\rightarrow A(X_{H})$ we denote the unique 
homomorphism whose pullback to $\pi_{1}(M)$ coincides with $\tau$.
Consider the action of
$H$ on $V=\sqcup_{i\in I}U_{i}\times \pi_{1}(M)$ induced by the left 
action of $H$ on $\pi_{1}(M)$ and let
$V_{H}=\sqcup_{i\in I}U_{i}\times X_{H}$ be the corresponding quotient set.
Define the equivalence relation
$U_{i}\times X_{H}\ni x\times \tau_{Q(H)}(\widetilde c_{ij})(h)
\sim x\times h\in U_{j}\times X_{H}$ 
with the same $\{c_{ij}\}$ as in the definition of $M(e)$.
The corresponding quotient space is a holomorphic bundle with fibre $X_{H}$ 
biholomorphic to $M(H)$. \\
{\bf (b)}\ We retain the notation of example (a). Let 
$B$ be a complex Banach space and
$GL(B)$ be the group of invertible
bounded linear operators $B\to B$. Consider a homomorphism 
$\rho: G\rightarrow GL(B)$. Without loss of generality we assume that 
$Ker(\rho)=e$, for otherwise we can pass to the corresponding quotient
group. The {\em holomorphic Banach vector bundle
$E_{\rho}\rightarrow M$ associated with $\rho$ } is defined as 
the quotient of $\sqcup_{i\in I} U_{i}\times B$ by the equivalence
relation
$U_{i}\times B\ni x\times\rho(c_{ij})(w)\sim x\times w\in U_{j}\times B$
for any $x\in U_{i}\cap U_{j}$. Let us illustrate this construction by 
an example. 

Let $\phi: X_{H}\to\Re^{+}$ ($X_{H}:=\pi_{1}(M)/H$) be a function satisfying
\begin{equation}\label{dilat}
\phi(\tau(h)(x))\leq c_{h}\phi(x)\ ,\ \ \ x\in X_{H}\ ,\ 
h\in\pi_{1}(M)\ ,
\end{equation}
where $c_{h}$ is a constant depending on $h$. By $l_{p,\phi}(X_{H})$, 
$1\leq p\leq\infty$, we denote
the Banach space of complex functions $f$ on $X_{H}$ with norm
\begin{equation}\label{p1}
||f||_{p,\phi}:=\left(\sum_{g\in X_{H}}|f(g)|^{p}\phi(g)\right)^{1/p} .
\end{equation}
Then according to (\ref{dilat}) the map $\rho$ defined by the
formula $[\rho(g)(f)](x):=f(\tau(g)(x))$, $g\in\pi_{1}(M)$, $x\in X_{H}$, is 
a homomorphism of $\pi_{1}(M)$ into $GL(l_{p,\phi}(X_{H}))$. By 
$E_{p,\phi}(X_{H})$ we denote the holomorphic
Banach vector bundle associated with this $\rho$.}
\end{E}
{\bf 3.2.}
We retain the notation of Example \ref{constbun}. Let $r:M'\to M$ be a 
covering where $M'=M(H)$ (i.e., $\pi_{1}(M')=H$).
Assume that $M$ satisfies condition (\ref{e1}), i.e., 
$M\subset\subset N$ and $\pi_{1}(M)=\pi_{1}(N)$.
Then there is an embedding $M(H)\hookrightarrow N(H)$.
(Without loss of generality we consider $M(H)$ as an
open subset of $N(H)$.) Let $\{V_{i}\}_{i\in I}$ be
a finite acyclic open cover of $\overline{M}$ by relatively compact sets.
We set $U_{i}:=V_{i}\cap M$ and consider the open cover 
${\cal U}=\{U_{i}\}_{i\in I}$ of $M$.
Then as in Example \ref{constbun} (a) we can define $M(H)$ by a cocycle
$c=\{c_{ij}\}\in Z^{1}_{{\cal O}}({\cal U},\pi_{1}(M))$.

Further, let $\psi:N(H)\to\Re_{+}$ be a function such that $\log\psi$ is 
uniformly continuous with respect to the path metric induced by a Riemannian 
metric pulled back from $N$. Fix a point $z_{0}\in M$ and identify 
$r^{-1}(z_{0})$ with $z_{0}\times X_{H}$\ ($X_{H}:=\pi_{1}(M)/H$). We 
define the function $\phi: X_{H}\to\Re_{+}$ by the 
formula
$$
\phi(x):=\psi(z_{0},x)\ ,\ \ \ x\in X_{H}\ .
$$
It was proved in [Br1, Lemma 2.3] that
{\em $\phi$ satisfies inequality} (\ref{dilat}). Then
the bundle $E_{p,\phi}(X_{H})$ is well defined.
By definition, any holomorphic section of this bundle is determined 
by a family $\{f_{i}(z,g)\}_{i\in I}$ of holomorphic functions on
$U_{i}$ with values in $l_{p,\phi}(X_{H})$ satisfying
$$
f_{i}(z,\tau(c_{ij})(h))=f_{j}(z, h)\ \ \ {\rm for\ any}\ \ \ 
z\in U_{i}\cap U_{j}\ .
$$
We introduce the Banach space $B_{p,\phi}(X_{H})$ of {\em
bounded holomorphic sections}
$f=\{f_{i}\}_{i\in I}$ of $E_{p,\phi}(X_{H})$ with norm
\begin{equation}\label{se1}
|f|_{p,\phi}:=\sup_{i\in I, z\in U_{i}}||f_{i}(z,\cdot)||_{p,\phi}\ .
\end{equation}
(Here $||\cdot ||_{p,\phi}$ is the norm on $l_{p,\phi}(X_{H})$, see 
(\ref{p1}).)

Further, let $f\in {\cal H}_{p,\psi}(M(H))$ (see section 2.2 for the
definition). We define the family
$\{f_{i}\}_{i\in I}$ of functions on $U_{i}$ with values in the space of 
functions on $X_{H}$ by the formula
\begin{equation}\label{cor1}
f_{i}(z,g):=f(z,g)\ ,\ \ \ z\in U_{i}\ ,\ i\in I\ ,\ g\in X_{H}\ .
\end{equation}
It was established in [Br1, Proposition 2.4] that {\em the
correspondence $f\mapsto \{f_{i}\}_{i\in I}$ determines an isomorphism
of Banach spaces} $D:{\cal H}_{p,\psi}(M(H))\to B_{p,\phi}(X_{H})$.
(Here $D$ is an isometry for $\psi\equiv 1$.) 

Next, suppose that $\{x_{n}\}_{n\geq 1}\subset M$ converges to $x\in M$. Then
for sufficiently big $n$ we can arrange $r^{-1}(x_{n})$ and $r^{-1}(x)$ in 
sequences 
$\{y_{in}\}_{i\geq 1}$ and $\{y_{i}\}_{i\geq 1}$ such that every
$\{y_{in}\}$ converges to $y_{i}$ as $n\to\infty$. For such $n$ we
define maps $\tau_{n}(x): r^{-1}(x)\to r^{-1}(x_{n})$ so that 
$\tau_{n}(y_{i})=y_{in}$, $i\in\N$. Below, $\tau_{n}^{*}$ denotes
the transpose map generated by $\tau_{n}$ on functions defined on
$r^{-1}(x_{n})$ and $r^{-1}(x)$.
\begin{D}\label{cont1} 
Let $X\subset M$ be a subset. We say that a function $f$ on $r^{-1}(X)$
belongs to the class $C_{p,\psi}(r^{-1}(X))$ if
\begin{itemize}
\item[{\rm (1)}]
$f|_{r^{-1}(x)}\in l_{p,\psi,x}(M')$ for all $x\in X$ and 
\item[{\rm (2)}]
for any $x\in X$ and 
any sequence $\{x_{n}\}\subset X$ converging to $x$ the sequence of functions
$\{\tau_{n}^{*}(f|_{r^{-1}(x_{n})})\}$ converges to 
$f|_{r^{-1}(x)}$ in the norm of $l_{p,\psi,x}(M')$.
\end{itemize}
\end{D}

By $C_{p,\psi}^{b}(r^{-1}(X))$ we denote the Banach space of functions
$f\in C_{p,\psi}(r^{-1}(X))$ with norm
\begin{equation}\label{norm}
|f|_{p,\psi}^{X}:=\sup_{x\in X}|f|_{r^{-1}(x)}|_{p,\psi,x}\ .
\end{equation}
Note that if $X\subset M$ is compact, then  $|f|_{p,\psi}^{X}<\infty$ for
every $f\in C_{p,\psi}(r^{-1}(X))$.

Comparing with the above definition of $D$ one determines a similar map for 
$C_{p,\psi}^{b}(r^{-1}(X))$. This gives
an isomorphism $D: C_{p,\psi}^{b}(r^{-1}(X))\to CB_{p,\phi}^{X}(X_{H})$ where
$CB_{p,\phi}^{X}(X_{H})$ is the Banach space of 
{\em bounded continuous sections}
of $E_{p,\phi}(X_{H})|_{X}$ with norm defined as in (\ref{se1}).\\
{\bf 3.3.} Most of our proofs are based on Theorem 1.3 of [Br1]. In its
proof we use the above isomorphisms $D$ and Cartan's A and B theorems
for coherent Banach vector sheaves (see [B]). Let us formulate this
result.

Suppose that $r: M'\to M$ is a covering with $M$ satisfying (\ref{e1}).
We define $\psi: M'\to\Re_{+}$ as in section 3.2. Also, we define
${\cal H}_{p,\psi}(M')$ and $l_{p,\psi,x}(M')$ as in section 2.3.
For Banach spaces $E$ and $F$ by ${\cal B}(E,F)$ we denote
the space of all linear bounded operators $E\to F$ with norm $||\cdot||$.
\begin{Th}\label{previous}
For any $p\in [1,\infty]$ there is a family 
$\{L_{z}\in {\cal B}(l_{p,\psi,z}(M'),
{\cal H}_{p,\psi}(M'))\}_{z\in M}$ holomorphic 
in $z$ such that 
$$
(L_{z}h)(x)=h(x)\ \ \ {\rm for\ any}\ \ \ h\in l_{p,\psi,z}(M')\ \ \
{\rm and}\ \ \ x\in r^{-1}(z)\ .
$$
Moreover,
$$
\sup_{z\in M}||L_{z}||<\infty\ .
$$
\end{Th}

The following facts are simple corollaries of this result.

Suppose that $X\subset M$ and $f\in C_{p,\psi}^{b}(r^{-1}(X))$. We define
the function $F$ on $X\times M'$ by the formula
\begin{equation}\label{newfun}
F(x,z):=(L_{x}(f|_{r^{-1}(x)}))(z)\ ,\ \ \ \ \ x\times z\in X\times M'\ .
\end{equation}
Then $F$ is continuous and $F(x,\cdot)\in {\cal H}_{p,\psi}(M')$ for every
$x$. Moreover, if $X$ is open and 
$f\in {\cal H}_{p,\psi}(r^{-1}(X))$, then $F\in {\cal O}(X\times M')$ and the
map $X\to {\cal H}_{p,\psi}(M')$,
$x\mapsto F(x,\cdot)$, is holomorphic.

We can also express $F$ in local coordinates. Namely, take
$x\in X$ and let $U\subset M$ be a neighbourhood of $x$ 
biholomorphic to an open Euclidean ball. Then 
$r^{-1}(U)=\sqcup_{y\in r^{-1}(x)}V_{y}$ and there are biholomorphisms
$s_{y}:U\to V_{y}$ such that $r\circ s_{y}=id$. Now, the restriction of
$f\in C_{p,\psi}^{b}(r^{-1}(X))$ to $r^{-1}(U\cap X)$ can be written as
\begin{equation}\label{new}
f(z)=\sum_{y\in r^{-1}(x)}f(z)\chi_{y}(z)\ ,\ \ \ z\in r^{-1}(U\cap X)\ ,
\end{equation}
where $\chi_{y}$ is the characteristic function of $V_{y}$. Let us
introduce the functions $\widetilde f_{y}$, $y\in r^{-1}(x)$, by the formulas
$$
\widetilde f_{y}(v)=f(s_{y}(v))\ ,\ \ \ v\in X\cap U\ .
$$
Then we have
\begin{equation}\label{new1}
f(z)=\sum_{y\in r^{-1}(x)}\widetilde f_{y}(v)\chi_{y}(z)\ ,\ \ \ \
v=r(z)\in U\cap X\ .
\end{equation}
Consider the series 
\begin{equation}\label{series}
\sum_{y\in r^{-1}(x)}
\widetilde f_{y}(v)L_{v}(\chi_{y}|_{r^{-1}(v)})\ ,\ \ \ v\in  U\cap X\ ,
\end{equation}
with $L_{v}$ as in Theorem \ref{previous}.
\begin{Proposition}\label{conv}
For $p\in [1,\infty)$ the series in (\ref{series}) converges
in ${\cal H}_{p,\psi}(M')$ to $F(v,\cdot):=L_{v}(f|_{r^{-1}(v)})$
uniformly on every compact subset of $U\cap X$. If $p=\infty$ and 
$f\in C_{1,1}^{b}(r^{-1}(X))$ then 
this series also converges in ${\cal H}_{\infty,\psi}(M')$ to
$F(v,\cdot)$ uniformly on every compact subset of $U\cap X$.
\end{Proposition}
{\bf Proof.} Suppose that $p\in [1,\infty)$ and 
$f\in C_{p,\psi}^{b}(r^{-1}(X))$.
Let $C\subset U\cap X$ be a compact subset. By the definition the function
$\Phi:U\cap X\to l_{p,\psi}(X_{H})$, $z\mapsto\widetilde f_{\cdot}(z)$, is
continuous (here we identify $r^{-1}(x)$ with $X_{H}$).
Thus $\Phi(C)\subset l_{p,\psi}(X_{H})$ is compact. Fix a 
family $\{X_{i}\}_{i\in\N}$ of finite subsets of $X_{H}$ such that
$X_{i}\subset X_{i+1}$ for any $i$ and $\cup_{i=1}^{\infty}X_{i}=X_{H}$.
Let $V_{i}\subset l_{p,\psi}(X_{H})$ be a finite-dimensional subspace
generated by functions $\delta_{z}$ on $X_{H}$ with $z\in X_{i}$. Here
$\delta_{z}(v)=1$ if $v=z$ and $\delta_{z}(v)=0$ if $v\neq z$. Then
$\cup_{i=1}^{\infty}V_{i}$ is everywhere dense in $l_{p,\psi}(X_{H})$
(since $1\leq p<\infty$). This and compactness of $\Phi(C)$ imply
that for any $\epsilon>0$ there exists an integer $l$ such that
$\Phi(C)\subset V_{l}+B_{\epsilon}$ where $B_{\epsilon}$ is the open ball
in $l_{p,\psi}(X_{H})$ centered at $0$ of radius $\epsilon$. By
$p_{l}: l_{p,\psi}(X_{H})\to V_{l}$ we denote the projection sending
$v=\sum_{x\in X_{H}}v_{x}\delta_{x}\in l_{p,\psi}(X_{H})$ to
$\sum_{x\in X_{l}}v_{x}\delta_{x}\in V_{l}$ (here all $v_{x}\in\Co$).
Then for $\epsilon$ as above and every $v\in\Phi(C)$ we have 
$||v-p_{l}(v)||_{p,\phi}<\epsilon$. From this by (\ref{new1}), identifying 
$r^{-1}(x)$ with $X_{H}$, we obtain
\begin{equation}\label{new2}
\sup_{v\in C}||f|_{r^{-1}(v)}-\sum_{y\in X_{l}}\widetilde f_{y}(v)
\chi_{y}|_{r^{-1}(v)}||_{p,\phi}<\epsilon\ .
\end{equation}
Thus by the definition of operators $L_{v}$ (see Theorem \ref{previous})
\begin{equation}\label{new3}
\sup_{v\in C}|F(v,\cdot)-\sum_{y\in X_{l}}\widetilde f_{y}(v)
L_{v}(\chi_{y}|_{r^{-1}(v)})|_{p,\psi}^{M}< C\epsilon
\end{equation}
for some constant $C$. This implies the required uniform convergence for
$p\in [1,\infty)$.

For $p=\infty$ and 
$f\in C_{1,1}^{b}(r^{-1}(X))$
we obtain anew that $\Phi(C)\subset l_{1,1}(X_{H})$ is compact. Then
in the above notation we easily get $||v-p_{l}(v)||_{\infty,\phi}<\epsilon$ 
for any $v\in\Phi(C)$ (because $||\cdot||_{\infty,\phi}=||\cdot||_{\infty,1}
\leq ||\cdot||_{1,1}$). Thus (\ref{new2}) is also valid for
$p=\infty$. This gives (\ref{new3}) with $p=\infty$.\ \ \ \ \ $\Box$
\sect{\hspace*{-1em}. Proofs of Theorem \ref{te1} and Corollary \ref{co1}.}
{\bf Proof of Theorem \ref{te1}.}
Let $M\subset\subset\widetilde M\subset N$ be complex manifolds such that
$\pi_{1}(M)=\pi_{1}(N)$ and $\widetilde M$ is Stein. Consider
an unbranched covering $r:N'\to N$ of $N$
and the corresponding coverings $M'=r^{-1}(M)$ and 
$\widetilde M'=r^{-1}(\widetilde M)$ of $M$ and $\widetilde M$. 
According to Example \ref{constbun} (a) $\widetilde M'$ is defined
on an open cover ${\cal U}=\{U_{i}\}_{i\in I}$ of $\widetilde M$ by
sets biholomorphic to open Euclidean balls by a locally constant cocycle
$\{\widetilde c_{ij}\}\in Z_{{\cal O}}^{1}({\cal U},Q(H))$. 
(Here we retain the notation of Example \ref{constbun} (a) so that
$\widetilde M'=\widetilde M(H)$.)
Using this construction we 
identify $r^{-1}(U_{i})$ with $U_{i}\times X_{H}$ ($X_{H}=\pi_{1}(M)/H$). 
Also, we choose some points
$z_{i}\in U_{i}$ and assume that diameters of all $U_{i}$ in the path metric 
on $N$ induced by a Riemannian metric are uniformly bounded by a constant.

Let $\phi:N'\to\Re$ be a function uniformly continuous with respect to
the path metric induced by the Riemannian metric pulled back from $N$.
For every $i\in I$ we define a function $\phi_{i}:r^{-1}(U_{i})\to\Re$
by the formula
$$
\phi_{i}(z,g):=\phi(z_{i},g)\ ,\ \ \ z\times g\in r^{-1}(U_{i})\ .
$$
Then from the uniform continuity of $\phi$ and boundedness of $diam(U_{i})$
for all $i$ we obtain that there exists a constant $c$ such that
\begin{equation}\label{phi}
|\phi(v)-\phi_{i}(v)|\leq c\ \ \ {\rm for\ every}\ \ \ 
v\in r^{-1}(U_{i})\ ,\ i\in I\ .
\end{equation}
Define a locally constant cocycle $\phi_{ij}$ on the open
cover $\{r^{-1}(U_{i})\}$ of $\widetilde M'$ by the formula
$$
\phi_{ij}(v)=\phi_{i}(v)-\phi_{j}(v)\ \ \ {\rm for}\ \ \ 
v\in r^{-1}(U_{i}\cap U_{j})\ .
$$
Then from (\ref{phi}) by the triangle inequality we get 
\begin{equation}\label{phi1}
\sup_{i, j, v}|\phi_{ij}(v)|\leq 2c\ .
\end{equation}
This inequality implies that rewriting cocycle $\{\phi_{ij}\}$ in
the coordinates on $\widetilde M$ (i.e., taking its direct image 
$\{r_{*}(\phi_{ij})\}$ with respect to $r$) 
we can regard it as a holomorphic cocycle on the cover ${\cal U}$ with 
values in the Banach vector bundle $E_{\infty, 1}(X_{H})$ with fibre 
$l_{\infty,1}(X_{H})$ defined on $\widetilde M$
(see Example \ref{constbun} (b)). This correspondence is 
described in [Br2, Proposition 2.4]. Since ${\cal U}$ is acyclic and
$\widetilde M$ is Stein, from the above construction, a version of 
Cartan's B theorem for coherent Banach sheaves (see [B]), and the classical 
Leray theorem we obtain as in [Br2] that there are 
holomorphic functions $f_{i}\in {\cal O}(r^{-1}(U_{i}))$ such that
\begin{itemize}
\item[(1)]
for every compact set $K\subset U_{i}$,
$$
\sup_{y\in r^{-1}(K)}|f_{i}(y)|<\infty\ ;
$$
\item[(2)]
$$
f_{i}(z)-f_{j}(z)=\phi_{ij}(z)\ \ \ {\rm for}\ \ \ z\in 
r^{-1}(U_{i}\cap U_{j})\ .
$$
\end{itemize}
Let ${\cal V}=\{V_{j}\}_{j\in J}$ be a refinement of ${\cal U}$ such that 
every $V_{j}$ is open and relatively compact in some $U_{i_{j}}$. Then 
condition (1) implies that
\begin{equation}\label{bou2}
\sup_{y\in r^{-1}(V_{j})}|f_{i_{j}}(y)|<\infty\ .
\end{equation}
Finally, define a function $\widetilde f\in {\cal O}(\widetilde M')$ by 
the formula
\begin{equation}\label{func}
\widetilde f(z):=\phi_{i}(z)-f_{i}(z)\ ,\ \ \ z\in r^{-1}(U_{i})\ .
\end{equation}
Since $\overline{M}\subset\widetilde M$ is a compact set, there is
a finite subcover of ${\cal V}$ that covers $\overline{M}$. From here,
(\ref{bou2}) and (\ref{phi}) for the restriction 
$f_{\phi}:=\widetilde f|_{\overline{M'}}$ we obtain (for some $C$)
\begin{equation}\label{func1}
|f_{\phi}(z)-\phi(z)|<C\ \ \ {\rm and}\ \ \ |df_{\phi}(z)|<C
\ \ \ {\rm for\ any}\ \ \ z\in M'\ .
\end{equation}

The proof of the theorem is complete.\ \ \ \ \ $\Box$\\
{\bf Proof of Corollary \ref{co1}.} We retain the notation of the proof of
Theorem \ref{te1}. 

Let $d_{o}:=d(o,\cdot)$ be the distance on $N'$ 
from a fixed point $o\in M'$ and let $f:=f_{\phi}$ be the function from 
Theorem \ref{te1} for $\phi=d_{o}$. Consider a finite open cover 
$\{U_{i}\}_{i=1}^{l}$ of $\overline{M}$ such that every 
$U_{i}\subset\subset\widetilde M$ is 
biholomorphic to an open Euclidean ball. As above we identify
$r^{-1}(U_{i})\subset N'$ with $U_{i}\times X_{H}$. 
Fix an element $e\in X_{H}$ and set 
$o_{i}(z)=(z,e)\in r^{-1}(U_{i})$ for every $z\in U_{i}$, $1\leq i\leq l$,
and $d_{o_{i}(z)}(v):=d(o_{i}(z),v)$, $v\in N'$. 
Then from compactness of every $\overline{U_{i}}$ by the triangle inequality 
we get
\begin{equation}\label{comp1}
|d_{o_{i}(z)}(v)-d_{o}(v)|\leq a\ ,\ \ \ 1\leq i\leq l\ ,
\end{equation}
for some constant $a$. By $B_{z,i}(R)$ we denote the open ball on $r^{-1}(z)$
of radius $R$ centered at $o_{i}(z)$
with respect to the induced metric $d|_{r^{-1}(z)}$. Also, by $\# A$ we
denote the number of elements of $A$. Now we prove
\begin{Lm}\label{number}
There is $k\in\N$ such that 
$$
\# B_{z,i}(R)\leq e^{kR}\ ,\ \ \ 1\leq i\leq l\ .
$$
\end{Lm}
{\bf Proof.}
Let $\widetilde r:N_{u}\to N$ be the universal covering of $N$ and
$r':N_{u}\to N'$ be the intermediate covering, i.e., $\widetilde r=r\circ r'$.
We equip $N_{u}$ with the path metric $\widetilde d$ induced by the 
Riemannian metric pulled back from $N$, the same as in the definition of the 
metric $d$ on $N'$. Let $\widetilde o_{i}(z)\in \widetilde r^{-1}(z)$ be 
such that $r'(\widetilde o_{i}(z))=o_{i}(z)$. By 
$\widetilde B_{\widetilde z,i}(R)$ we denote the open ball on 
$\widetilde r^{-1}(z)$ of radius $R$ centered at 
$\widetilde o_{i}(z)$ with respect to the metric 
$\widetilde d|_{\widetilde r^{-1}(z)}$. 
Let us check that $r'(\widetilde B_{\widetilde z,i}(R))=B_{z,i}(R)$.

Indeed, let $y\in\widetilde B_{\widetilde z,i}(R)$ and $\gamma_{y}$ be a 
path joining $\widetilde o_{i}(z)$ and $y$ in $N_{u}$ whose length
is less than $R$ (such a path exists by the definition of $\widetilde d$).
Then $r'(\gamma_{y})$ is a path joining $o_{i}(z)$ and $r'(y)$ in $N'$.
By the definition of the metrics on $N_{u}$ and $N'$ the length
of $r'(\gamma_{y})$ does not exceed the length of $\gamma_{y}$.
In particular, it is less than $R$. Thus $d_{o_{i}(z)}(r'(y))<R$, i.e.,
$r'(y)\in B_{z,i}(R)$. Conversely, let $w\in B_{z,i}(R)$ and let $\gamma_{w}$
be a path in $N'$ joining $o_{i}(z)$ and $w$ with length less than $R$. By the
covering homotopy theorem (see, e.g., [Hu, Chapter III]) there is a 
path $\widetilde\gamma_{w}\subset N_{u}$ that covers $\gamma_{w}$ and
joins $\widetilde o_{i}(z)$ with some point $\widetilde w$ such that
$r'(\widetilde w)=w$. Moreover, by the definition, the length of 
$\widetilde\gamma_{w}$ is the same as the length of $\gamma_{w}$. In 
particular, it is less than $R$. Thus 
$\widetilde w\in\widetilde B_{\widetilde z,i}(R)$. This shows that
$r'(\widetilde B_{\widetilde z,i}(R))=B_{z,i}(R)$. In turn, the latter implies
that
\begin{equation}\label{compare}
\# B_{z,i}(R)\leq\#\widetilde B_{\widetilde z,i}(R)\ .
\end{equation}

Next, let $A$ be a finite set of generators of $\pi_{1}(N)$ (recall that
condition (\ref{e1}) implies that $\pi_{1}(N)$ is finitely generated).
By $d_{w}$ we denote the word metric on $\pi_{1}(N)$ with respect to $A$. 
Now, from compactness of every $\overline{U_{i}}$ by the 
\v{S}varc-Milnor lemma (see, e.g., [BH, p.140])
we obtain that there exists a constant $c$ such that for any
$z\in U_{i}$, $1\leq i\leq l$, and $g, h\in\pi_{1}(N)$,
\begin{equation}\label{comp2}
c^{-1}d_{w}(g,h)\leq \widetilde d((z,g),(z,h))\leq cd_{w}(g,h)
\end{equation}
(Here we identify $\widetilde r^{-1}(U_{i})$ with $U_{i}\times\pi_{1}(N)$
as in Example \ref{constbun} (a).) Let $B_{R}\subset\pi_{1}(N)$ be the
open ball of radius $R$ centered at $1$ with respect to $d_{w}$. Then
there is a natural number $\widetilde k$ such that 
\begin{equation}\label{comp2'}
\#B_{R}\leq e^{\widetilde k R}\ \ \ {\rm for\ any}\ \ \ R\geq 0\ .
\end{equation}
From here, (\ref{comp2}) and (\ref{compare}) we get for $k:=\widetilde k c$
$$
\#B_{z,i}(R)\leq e^{k R}\ ,\ \ \ 1\leq i\leq l\ . \ \ \ \ \ \Box
$$

We proceed with the proof of the corollary.
Let us  define $\alpha:=\frac{k+1}{2}$ and prove that 
$F=e^{-\alpha f}\in {\cal H}_{2,1}(M')\cap C(\overline{M'})$.
Since ${\cal H}_{2,1}(M')\hookrightarrow H_{1}^{2}(M')$ (see section 2.3),
this implies the required statement.

Let $z\in U_{i}$ for some $1\leq i\leq l$. We will estimate
$|F|_{2,1,z}$ (see (\ref{e5})). By the definition using inequalities
(\ref{func1}), (\ref{comp1}) and Lemma  \ref{number} we obtain
$$
\begin{array}{c}
\displaystyle
|F|_{2,1,z}^{2}=\sum_{y\in r^{-1}(z)}|e^{-\alpha f(y)}|^{2}\leq
e^{2\alpha (a+C)}\cdot\sum_{y\in r^{-1}(z)}e^{-2\alpha d_{o_{i}(z)}(y)}\leq\\
\\
\displaystyle
e^{2\alpha (a+C)}\cdot\sum_{R=0}^{\infty}e^{-2\alpha R}\cdot 
\# B_{z,i}(R)\leq
e^{2\alpha (a+C)}\cdot\sum_{R=0}^{\infty}e^{(-2\alpha +k)R}=
\frac{e^{2\alpha (a+C)+1}}{e-1}\ .
\end{array}
$$
Therefore
$$
|F|_{2,1}^{M}:=\sup_{z\in M}|F|_{2,1,z}\leq
\left(\frac{e^{2\alpha (a+C)+1}}{e-1}\right)^{1/2}\ .
$$
This shows that $F\in {\cal H}_{2,1}(M')\cap C(\overline{M'})$.\ \ \ \ \ 
$\Box$
\begin{R}\label{re1'}
{\rm (1) Using some construction from [Br2] one can prove that the 
constant $C$ in Theorem \ref{te1} for $\phi=d_{o}$ (see (\ref{func1})) can be 
chosen independent of the covering $r:M'\to M$. It depends only on 
$M$, $\widetilde M$ and the Riemannian metric on $N$.\\
(2) Consider the holomorphic map $f:M'\to\Co$ with $f$ as in Corollary 
\ref{co1}. Then
$$
f(M')\subset S:=\{z\in\Co\ :\ |Im\ \!z|<C\ ,\ -C< Re\ \!z<\infty\}
$$
where $C$ is the constant in Theorem \ref{te1} for $\phi=d_{o}$.
Let $B_{t}=\{x\in M'\ :\ d_{o}(x)<t\}$ be the open ball in $M'$ centered
at $o$ of radius $t$ and
$S_{R}:=\{z\in S\ :\ Re\ \! z\geq R\}$.
Then 
$$
f^{-1}(S_{R})\subset M'\setminus B_{R-C}\ \ \ {\rm for}\ \ \ R>C\ \ \
{\rm and}\ \ \ f^{-1}(S\setminus S_{R})\subset B_{R+C}\ .
$$
}
\end{R}
Using such $f$ one can construct holomorphic
functions on $M'$ decreasing faster than the function $F$ from
Corollary \ref{co1}. Actually, let $l:\Re_{+}\to\Re_{+}$ be a continuous
function monotonically increasing for $x\geq R_{0}$. Consider a holomorphic
function $g$ on $S$ satisfying
\begin{equation}\label{bound}
\log|g(x+iy)|\geq l(x)\ \ \ {\rm for}\ \ \ x\geq R_{0}\ \ \ \ {\rm and}
\ \ \ \ \inf_{z\in S}|g(z)|>0\ .
\end{equation}
Then one can easily check that the function
$G=g\circ f\in {\cal O}(M')\cap C(\overline{M'})$ satisfies
$$
|G(z)|\geq e^{l(d_{o}(z)-C)}\ \ \ {\rm for}\ \ \ d_{o}(z)\geq R_{0}+C\ .
$$
In particular, $H=1/G\in  
{\cal O}(M')\cap C(\overline{M'})$ 
satisfies (for some $c_{1}>0$)
\begin{equation}\label{bound1}
|H(z)|\leq c_{1}e^{-l(d_{o}(z)-C)}\ ,\ \ \
z\in\overline{M'}\setminus B_{C}\ .
\end{equation}
Observe that the Harnack inequality for positive harmonic functions implies
for $g$ as in (\ref{bound}) (for some positive $\widetilde c_{1}$,
$\widetilde c_{2}$)
$$
l(x)\leq \log|g(x)|\leq\widetilde c_{1}e^{\widetilde c_{2}x}\ ,\ \ 
\ x\geq R_{0}\ .
$$
This and the properties of $f$ impose the following restriction on the decay 
of $H$:
\begin{equation}\label{bound2}
|H(z)|\geq e^{-\widetilde c_{3}e^{\widetilde c_{2}d_{o}(z)}}\ ,\ \
\ z\in M'\ .
\end{equation}
\begin{E}\label{ex1}
{\rm As the function $g$ in (\ref{bound}) one can take, e.g.,
$g(z)=e^{z^{n}}$ for $n\in\N$ (in this case 
$l(x)=(1-\epsilon)x^{n}$ for any $\epsilon>0$), or
$g(z)=e^{C_{1}e^{C_{2}z}}$ for $C_{1}>0$ and $0<C_{2}<\frac{\pi}{2C}$ with
$C$ as above (in this case $l(x)=C_{1}\cos(C_{2}C)e^{C_{2}x}$).
For the latter example estimate (\ref{bound1})
shows that the lower bound (\ref{bound2}) of the decay of $H$
is attainable.}
\end{E}
\sect{\hspace*{-1em}. Proof of Theorem \ref{te2}.}
Suppose that $M\subset\subset N$ are domains in a Stein manifold, 
$M$ is strictly pseudoconvex and $\pi_{1}(M)=\pi_{1}(N)$. Using the
Remmert embedding theorem (see, e.g., [GR]) we may assume without loss of 
generality that $N$ is a domain in a closed complex submanifold of some 
$\Co^{k}$. Let $r:N'\to N$ be an unbranched covering of $N$. As usual,
we set $M'=r^{-1}(M)$ and
$bM'=r^{-1}(bM)$ where $bM$ is the boundary of $M$. 
We must show that every point in $bM'$ is a peak point for 
${\cal H}_{p,1}(M')$, $1\leq p\leq\infty$. In our proof we use a result on
uniform estimates for solutions of certain $\overline\partial$-equations
on $M'$. To its formulation we first introduce the corresponding class of
$(0,1)$-forms on $M'$.

Let $\{V_{i}\}_{i\in I}$ be a finite acyclic open cover of $\overline{M}$
by relatively compact complex coordinate systems. We set $U_{i}\cap M$
and consider the open cover ${\cal U}=\{U_{i}\}_{i\in I}$ of $M$. Let
$X_{H}$ be the fibre of $N'$ with $H=\pi_{1}(N')$.
Using the construction of Example \ref{constbun} (a)
we identify $r^{-1}(U_{i})$ with $U_{i}\times X_{H}$. 
Let $\omega$ be a $(0,1)$-form on $M'$. Then in local coordinates on 
$r^{-1}(U_{i})$ it is presented as 
$$
\omega(v,x)=\sum_{j=1}^{n}a_{j}(v,x)\ \!d\overline{v}_{j}\ \ \
{\rm for}\ \ \ v\times x\in U_{i}\times X_{H}\ 
$$
where $v=(v_{1},\dots, v_{n})$ are coordinates on $U_{i}$. 
Consider every $a_{j}$ as a function on $U_{i}$ with values in the
space of functions on $X_{H}$. We assume that for every $i\in I$
\begin{equation}\label{ass1}
a_{j}\in C^{\infty}(U_{i}, l_{p,1}(X_{H}))\
,\ 1\leq j\leq n\ .
\end{equation}
Then for such an $\omega$ its direct
image $r_{*}(\omega)$ is a bounded $C^{\infty}$ form with values in the
Banach vector bundle $E_{p,1}(X_{H})$ (see section 3.1).
Also, we assume that the norm of $\omega$ defined by the formula
\begin{equation}\label{norm1}
||\omega||:=\sup_{i\in I}\max_{1\leq j\leq n}|a_{j}|_{p,1}^{U_{i}}\ 
\end{equation}
is finite. (Recall that $|\cdot|_{p,1}^{U_{i}}$ are norms on
$C_{p,1}(r^{-1}(U_{i}))$, see  (\ref{norm}).)
\begin{Proposition}\label{destim}
There is a constant $C>0$ and for each $\overline\partial$-closed 
$(0,1)$-form $\omega$ satisfying (\ref{ass1}) there is a function
$f\in C^{\infty}(M')\cap C_{p,1}^{b}(M')$ such that 
$$
\overline\partial f=\omega\ \ \ {\rm and}\ \ \ |f|_{p,1}^{M}\leq 
C||\omega||\ .
$$
\end{Proposition}
{\bf Proof.} We apply the operators $L_{z}$ from
Theorem \ref{previous} to $\omega$. Namely, 
let us define a form $\widetilde\omega$ on $M$ by the formula
$$
\widetilde\omega(v,z):=
\sum_{j=1}^{n}(L_{v}a_{j}(v,\cdot))(z)\ \!d\overline{v}_{j}\ \ \ {\rm for}
\ \ \ v\times z\in U_{i}\times M'\ ,\ i\in I\ .
$$
It is readily seen that $\widetilde\omega$ is a bounded 
$\overline\partial$-closed
$C^{\infty}$ form on $M$ with values in ${\cal H}_{p,1}(M')$. We define
the norm of $\widetilde\omega$ by
\begin{equation}\label{norm2}
||\widetilde\omega||:=\sup_{i\in I, v\in U_{i}}
\max_{1\leq j\leq n}|L_{v}a_{j}(v,\cdot)|_{p,1}^{M}
\end{equation}
where $|\cdot|_{p,1}^{M}$ is norm on ${\cal H}_{p,1}(M')$. Then according to
Theorem \ref{previous} there is a constant $c$ (independent of
$\omega$) such that
\begin{equation}\label{ineq1}
||\widetilde\omega||\leq c||\omega||\ .
\end{equation}

Further, we use Lemma 1 from [He]. According to this lemma {\em there exist
a strictly pseudoconvex domain $W\subset\Co^{k}$ with $C^{2}$ boundary
such that $W\cap N=M$ and a holomorphic map $\pi$ from a neigbourhood
$U(\overline{W})$ of $\overline{W}$ onto $U(\overline{W})\cap N$ such
that $\pi(W)=M$ and $\pi|_{U(\overline{W})\cap N}$ is the identity map.}

Using this result we obtain that the pullback $\pi^{*}\widetilde\omega$ with 
respect to $\pi$ is a bounded
$\overline\partial$-closed $C^{\infty}$ form on $W$ with values in
${\cal H}_{p,1}(M')$. Moreover, there is a constant $c'$ (depending on
$\pi$ and $W$) such that
\begin{equation}\label{ineq2}
||\pi^{*}\widetilde\omega||\leq c'||\widetilde\omega||\ .
\end{equation}
Here for $\pi^{*}\widetilde\omega(w,\cdot)=
\sum_{j=1}^{k}\widetilde a_{j}(w,\cdot)d\overline{w}_{j}$,\ 
$w=(w_{1},\dots, w_{k})\in\Co^{k}$, we define
$$
||\pi^{*}\widetilde\omega||:=\sup_{w\in W}\max_{1\leq j\leq k}
|\widetilde a_{j}(w,\cdot)|_{p,1}^{M}\ .
$$ 

In [SH] uniform estimates for solutions of 
$\overline{\partial}$-equations on so-called pseudoconvex polyhedra were
obtained by means of global integral formulas. This class contains,
in particular, strictly pseudoconvex domains with $C^{2}$ boundaries.
Note that the estimates in [SH] remain valid if one solves
Banach-valued $\overline{\partial}$-equations. Therefore from
the results of [SH] we obtain that {\em there exists a bounded
$C^{\infty}$ function $h$ on $W$ with values in ${\cal H}_{p,1}(M')$ such that
$\overline\partial h=\pi^{*}\widetilde\omega$. Moreover,
\begin{equation}\label{ineq3}
||h||\leq c''||\pi^{*}\widetilde\omega||
\end{equation}
for some $c''$ (depending on $W$ only). Here}
$$
||h||:=\sup_{w\in W}|h(w,\cdot)|_{p,1}^{M}\ .
$$

Finally, define a function $f$ on $M'$  by the formula
$$
f(z):=h(r(z),z)\ ,\ \ \ z\in M'\ .
$$
Using that $r$ is holomorphic, 
$\widetilde\omega(r(z),z)=\omega(z)$, $z\in M'$, and 
$(\pi^{*}\widetilde\omega)|_{M}=\widetilde\omega$ we easily conclude that
$\overline\partial f=\omega$. By the definition 
$f\in C^{\infty}(M')\cap C_{p,1}^{b}(M')$
and from (\ref{ineq1})-(\ref{ineq3}) we have (for some $C$)
$$ 
|f|_{p,1}^{M}\leq C||\omega||\ .\ \ \ \ \ \Box
$$
\begin{R}\label{high}
{\rm (1) An analogous to Proposition \ref{ass1} statement is valid for a
similar class of bounded
$\overline\partial$-closed $(0,q)$-forms on $M'$.\\
(2) Using the main result of [He] and the estimates from [SH] one
can show that the result of Proposition \ref{ass1} is valid also for 
coverings of non-degenerate pseudoconvex polyhedrons on Stein manifolds 
(see [He] and [SH] for the definition).}
\end{R}

We pass to the proof of Theorem \ref{te2}. 
Take a point $z\in bM'$ and set $v=r(z)\in bM$. Let $U\subset\subset N$ be a 
simply connected coordinate neighbourhood of $v$ and let $W\subset N'$
be the neighbourhood of $z$ such that $r: W\to U$ is biholomorphic.
Since $M$ is strictly pseudoconvex, $v$ is a peak point for 
${\cal O}(U\cap M)$ for a sufficiently small $U$. Moreover, for such $U$
we can find $\widetilde f\in {\cal O}(U\cap M)$  
with a peak point at $v$ such that $\widetilde f\in L^{q}(U\cap M)$ for all
$1\leq q<\infty$ (see [GHS, p.575]). 
Then $f:=(r^{*}\widetilde f)|_{W\cap M'}$ has a peak point at $z$ and
$f\in L^{q}(W\cap M')$ for all $1\leq q<\infty$. 
Next, let $\widetilde\rho\in C^{\infty}(U)$ be a cut-off 
function which equals 1 in a neighbourhood $O\subset\subset U$ of $v$ 
and 0 outside $U$. Consider its pullback 
$\rho:=(r^{*}\widetilde\rho)|_{W}\in C^{\infty}(W)$. Clearly the
$(0,1)$-form $\omega=\overline\partial(\rho f)$ on $M'$ satisfies conditions
of Proposition \ref{destim}.
Then this proposition implies that there exists a 
function $h\in C^{\infty}(M')\cap C_{p,1}^{b}(M')$ such that 
$\overline\partial h=\omega$.
Finally, consider the function $h_{z}:=\rho f-h$. Then 
$h_{z}$ is holomorphic,  has a peak point at $z_{0}$ and 
belongs to $H_{1}^{p}(M')$ for $1\leq p<\infty$ by the choice of $f$.
Also, for $p=\infty$ the function $h_{z}$ is bounded outside $W$.

The proof of the theorem is complete.\ \ \ \ \ $\Box$
\sect{\hspace*{-1em}. Proofs of Theorem \ref{te4} and Corollary \ref{co2}.}
{\bf Proof of Theorem \ref{te4}.}
Let $M\subset\subset\widetilde M\subset N$ be manifolds satisfying condition
(\ref{e1}) with $dim_{\Co}M\geq 2$. Let
$D\subset\subset M$ be an open subset whose boundary $bD$ is a connected
$C^{k}$ submanifold of $M$ $(1\leq k\leq\infty)$.
Consider a covering $r:N'\to N$ and
set $M'=r^{-1}(M)$, $D'=r^{-1}(D)$ and $bD'=r^{-1}(bD)$. Let
$\psi:N'\to\Re_{+}$ be such that $\log\psi$ is uniformly continuous with
respect to the path metric induced by a Riemannian metric pulled back from
$N$. Let $f\in C^{s}(bD')$, $0\leq s\leq k$, be a $CR$-function satisfying 
the hypotheses of Theorem \ref{te4}.\\
(A)  First, we will prove the theorem for $s=0$ under the 
additional assumption 
\begin{equation}\label{assum1}
f\in C_{p,\psi}(bD')\ \ \ {\rm for}\ \ \ p\in [1,\infty)\ \ \
{\rm and}\ \ \ f\in C_{1,1}(bD')\ \ \
{\rm for}\ \ \ p=\infty\ .
\end{equation}
(We use here that $C_{1,1}(bD')\subset C_{\infty,1}(bD')=
C_{\infty,\psi}(bD')$, see Definition \ref{cont1}.)

For a $CR$-function $f$ satisfying (\ref{assum1}) we define a continuous
${\cal H}_{p,\psi}(M')$-valued function $F$ on
$bD'$ by the formula
$$
F(v):=L_{v}(f|_{r^{-1}(v)})\ ,\ \ \ v\in bD\ ,
$$
where $L_{v}$ are operators from Theorem \ref{previous}, see section 3.3.
\begin{Lm}\label{cr}
$F$ is a ${\cal H}_{p,\psi}(M')$-valued continuous $CR$-function.
\end{Lm}
{\bf Proof.}
Let $U\subset M$ be a simply connected coordinate neighbourhood of a point 
$x\in bD$. It suffices to check that $F|_{U\cap bD}$ satisfies the required 
property. Note that by Proposition \ref{conv}
\begin{equation}\label{series1}
[F(v)](z)=\sum_{y\in r^{-1}(x)}\widetilde f_{y}(v)H_{y}(v,z)\ ,\ \ \ 
v\times z\in (U\cap bD)\times M'\ ,
\end{equation}
where $\widetilde f_{y}(v)=f(s_{y}(v))$, $v\in U\cap bD$, and
$s_{y}:U\to V_{y}$ is a biholomorphic map onto the connected 
component $V_{y}$ of $r^{-1}(U)$ containing $y$. By the 
definition of operators $L_{v}$ functions $H_{y}$ are restrictions to 
$bD\times M'$ of some holomorphic functions on $U\times M'$.
Moreover, Proposition \ref{conv}
implies that the series in (\ref{series1}) converges uniformly to
$F$ on every compact subset of $(U\cap bD)\times M'$. Next,
since $f|_{V_{y}\cap bD'}$ is a continuous $CR$-function, $\widetilde f_{y}$ 
is a 
continuous $CR$-function on $U\cap bD$. Also, by the definition of $H_{y}$, 
for a fixed $z\in M'$ every $H_{y}(\cdot, z)$ is a continuous $CR$-function 
on $U\cap bD$. Hence, $\widetilde f_{y}\cdot H_{y}(\cdot, z)$ is a 
continuous  $CR$-function on $U\cap bD$, as well. 
Indeed, for each  $(n,n-2)$-form 
$\omega$ with a compact support in $U$ we have
$$
\int_{U\cap bD}\widetilde f_{y}(v)\cdot H_{y}(v,z)\ \!
\overline\partial\omega(v)=\int_{U\cap bD}\widetilde f_{y}(v)\ \!
\overline\partial(H_{y}(v,z)\cdot\omega(v))=0
$$
because $\widetilde f_{y}$ is $CR$. Since the series in (\ref{series1})
converges uniformly to $[F(\cdot)](z)$ on every compact subset of
$(U\cap bD)\times z$, every $[F(\cdot)](z)$, $z\in M'$, is a continuous
$CR$-function on $U\cap bD$. This implies the required statement.
\ \ \ \ \ $\Box$

Further, since $[F(v)](z)$ from Lemma \ref{cr} is holomorphic in $z\in M'$, 
we can expand it in the Taylor series in a 
complex coordinate neighbourhood $U_{z}$, 
\begin{equation}\label{taylor}
[F(v)](w)=\sum_{0\leq |\alpha|<\infty}F_{\alpha}(v)w^{\alpha}\ ,\ \ \
v\in bD\ .
\end{equation}
Here $\alpha=(\alpha_{1},\dots,\alpha_{s})\in(\Z_{+})^{s}$,
$|\alpha|=\sum_{i=1}^{s}\alpha_{i}$,
$w^{\alpha}=w_{1}^{\alpha_{1}}\dots w_{n}^{\alpha_{n}}$ and
$w=(w_{1},\dots, w_{n})$ are coordinates on $U_{z}$ such that $w(z)=0$.
Now, from Lemma \ref{cr} follows that each $F_{\alpha}$ in (\ref{taylor})
is a continuous $CR$-function on
$bD$. Then by Theorem 3.14 of Harvey [Ha] for every 
$F_{\alpha}$ there exists a function
$\widetilde F_{\alpha}\in {\cal O}(D)\cap C(\overline{D})$ such that 
$\widetilde F_{\alpha}|_{bD}=F_{\alpha}$. Also, for a sufficiently small
$U_{z}$ using estimates of the Cauchy integrals for derivatives of a
holomorphic function and compactness of $bD$ we get from (\ref{taylor})
$$
M:=\sup_{\alpha,v\in bD}|F_{\alpha}(v)|<\infty\ .
$$
Thus by the maximum modulus principle 
$$
\sup_{\alpha,v\in\overline{D}}|\widetilde F_{\alpha}(v)|=M<\infty\ .
$$
The latter implies that for a sufficiently small $U_{z}$ the series
$$
\widetilde F_{z}(v,w)=
\sum_{0\leq |\alpha|<\infty}\widetilde F_{\alpha}(v)w^{\alpha}\ ,\ \ \
v\times w\in\overline{D}\times U_{z}\ ,
$$
converges absolutely and uniformly. Hence,
$\widetilde F_{z}\in {\cal O}(D\times U_{z})\cap C(\overline{D}\times U_{z})$.
Further, assume that for $y,z\in M'$ we have $U_{y}\cap U_{z}\neq\emptyset$.
Then for every $w\in U_{y}\cap U_{z}$ and $v\in bD$ 
$$
\widetilde F_{y}(v,w)-\widetilde F_{z}(v,w)=[F(v)](w)-[F(v)](w)=0\ .
$$
This leads to the identity
$$
\widetilde F_{y}(\cdot,w)=\widetilde F_{z}(\cdot,w)\ ,\ \ \ 
w\in U_{y}\cap U_{z}\ .
$$
Thus we can define a function $F\in {\cal O}(D\times M')\cap 
C(\overline{D}\times M')$ by the formula
\begin{equation}\label{tifunc}
\widetilde F(v,w):=\widetilde F_{z}(v,w)\ ,\ \ \ v\times w\in
\overline{D}\times U_{z}\ .
\end{equation}
\begin{Lm}\label{bound3}
$\widetilde F(v,\cdot)\in {\cal H}_{p,\psi}(M')$ for any $v\in D$.
\end{Lm}
{\bf Proof.}
Observe that the evaluation at $v\in D$ is a linear continuous functional
on the Banach space ${\cal O}(D)\cap C(\overline{D})$
equipped with supremum norm. Identifying
${\cal O}(D)\cap C(\overline{D})$ with its trace space on $bD$ and using
the Hahn-Banach and F. Riesz theorems we have
$$
h(v)=\int_{bD}h(\xi)d\mu(\xi)\ ,\ \ \ h\in {\cal O}(D)\cap C(\overline{D})\ ,
$$
where $\mu$ is a complex regular Borel measure on $bD$ with the total
variation $Var\ \!\mu=1$. Thus for every fixed $w\in M'$ we have
$$
\widetilde F(v,w)=\int_{bD}[F(\xi)](w)d\mu(\xi)\ .
$$
Now, by the definition of the norm on ${\cal H}_{p,\psi}(M')$
using the triangle inequality, the identity $\widetilde F(v,\cdot)=
F(v)$, $v\in bD$, and the fact that $F$ is a continuous
${\cal H}_{p,\psi}(M')$-valued function on $bD$ we obtain
$$
\begin{array}{c}
\displaystyle
|\widetilde F(v,\cdot)|_{p,\psi}^{M}:=\sup_{z\in M}
\left(\sum_{y\in r^{-1}(z)}|\widetilde F(v,y)|^{p}\psi(y)\right)^{1/p}\leq\\
\\
\displaystyle
\sup_{z\in M}\left(\sum_{y\in r^{-1}(z)}\left(\int_{bD}|\widetilde F(\xi,y)| 
|d\mu(\xi)|\right)^{p}\psi(y)\right)^{1/p}\leq\\
\\
\displaystyle
\sup_{z\in M}\left(\int_{bD}\left(\sum_{y\in r^{-1}(z)}|[F(\xi)](y)|^{p}
\psi(y)
\right)^{1/p}|d\mu(\xi)|\right)\leq \sup_{\xi\in bD}
|F(\xi)|_{p,\psi}^{M}
<\infty\ .\ \ \ \ \ \Box
\end{array}
$$

Further, set
\begin{equation}\label{recfunc}
f'(z):=\widetilde F(r(z),z)\ ,\ \ \ z\in \overline{D'}\ .
\end{equation}
Then using the inequalities of Lemma \ref{bound3} we get 
$$
\begin{array}{c}
\displaystyle
f'\in {\cal O}(D')\cap C(\overline{D'})\ ,\ \ \
f'|_{bD'}=F|_{bD'}=f\ , \ \ \ {\rm and}\\
\\
\displaystyle
|f'|_{p,\psi,z}:=\left(\sum_{y\in r^{-1}(z)}|f'(y)|^{p}\psi(y)\right)^{1/p}
\leq \sup_{\xi\in bD}|F(\xi)|_{p,\psi}^{M}\leq 
c|f|_{p,\psi}^{bD}\ ,\ \ \ z\in D\ , 
\end{array}
$$
see (\ref{norm}) for the definition of $|\cdot|_{p,\psi}^{bD}$. 
(Here the last inequality follows directly from Theorem \ref{previous}.)
The latter implies that 
$f'|_{r^{-1}(z)}\in l_{p,\psi,z}(M')$, $z\in D$, see section 2.3. Thus 
$f'\in {\cal H}_{p,\psi}(D')\cap C(\overline{D'})$ and
\begin{equation}\label{norm3}
|f'|_{p,\psi}^{D}:=\sup_{z\in D}|f'|_{p,\psi,z}\leq 
c\sup_{z\in bD}|f|_{p,\psi,z}\ (:=c|f|_{p,\psi}^{bD})\ .
\end{equation}

This completes the proof of the theorem for $s=0$ under  assumption 
(\ref{assum1}).\\
(B) Let us consider the general case of a continuous $CR$-function 
$f$ on $bD'$ satisfying 
\begin{equation}\label{general}
f|_{r^{-1}(x)}\in l_{p,\psi,x}(M')\ \ \ {\rm for\ any}\ \ \ x\in D \ \ \
{\rm and}\ \ \
m:=\sup_{x\in bD}|f|_{p,\psi,x}<\infty\ .
\end{equation}
According to Remark \ref{re1'} (2) and Example \ref{ex1} there is a constant 
$c>0$ such that for any $0<c_{1}<c$ and  $c_{2}>0$ there exists a function 
$F_{c_{1},c_{2}}\in {\cal O}(M')\cap C(\overline{M'})$ satisfying
\begin{equation}\label{double}
e^{-c_{3}e^{c_{1}d_{o}(z)}}\leq |F_{c_{1},c_{2}}(z)|\leq 
e^{-c_{2}e^{c_{1}d_{o}(z)}}\ \ \ {\rm for\ all}\ \ \ z\in \overline{M'}\ 
\end{equation}
with $c_{3}$ depending on $c_{2}$, $c_{1}$, $c$, $M$, $M'$ such that
$c_{3}\to 0$ as $c_{2}\to 0$. 
(Recall that $d_{o}$, $o\in M'$, is the 
distance on $N'$ defined as in Corollary \ref{co1}.)
Define a continuous $CR$-function  $f_{c_{1},c_{2}}$ by the formula
$$
f_{c_{1},c_{2}}(z):=f(z)F_{c_{1},c_{2}}(z)\ ,\ \ \ z\in bD'\ .
$$
\begin{Lm}\label{smooth}
$f_{c_{1},c_{2}}$ satisfies assumption (\ref{assum1}).
\end{Lm}
{\bf Proof.} Note that for any $l\in\N$ there is a nonnegative $r$ such 
that 
$$
e^{-c_{2}e^{c_{1}d_{o}(z)}}<e^{-l d_{o}(z)}\ \ \ {\rm for}\ \ \ 
d_{o}(z)>r\ ,\ z\in\overline{M'}\ .
$$ 
Let $U\subset\subset M$ be a neighbourhood of $\overline{D}$ and
$U'=r^{-1}(U)\subset M'$. From the above inequality
arguing as in the proof of Corollary \ref{co1} we obtain that 
$F_{c_{1},c_{2}}\in {\cal H}_{p,1}(U')$ for any $p\in [1,\infty]$.
Then from [Br, Proposition 2.4] follows that 
$F_{c_{1},c_{2}}|_{bD'}$  belongs to $C_{p,1}(bD')$ for all $p$.

Next, take a point $x\in bD$ and prove that $f_{c_{1},c_{1}}$ is
$C_{p,\psi}$-continuous over $x$. Let $U_{x}\subset M$ be a complex 
(simply connected) coordinate neighbourhood of $x$. We will  identify
$r^{-1}(U_{x})$ with $U_{x}\times X_{H}$ where $X_{H}$ is the fibre of
$r:M'\to M$. Consider a sequence $\{x_{n}\}\subset U_{x}\cap bD$ convergent to
$x$. For $g\in X_{H}$ put $a_{n}(g):=f(x_{n},g)$, 
$a(g):=f(x,g)$,
$b_{n}(g):=F_{c_{1},c_{2}}(x_{n},g)$, $b(g):=F_{c_{1},c_{2}}(x,g)$ and 
$c_{n}:=a_{n}b_{n}$, $c:=ab$.
Then we must check that
$$
\lim_{n\to\infty}|c-c_{n}|_{p,\psi,x}:=
\lim_{n\to\infty}\left(\sum_{g\in X_{H}}|c(g)-c_{n}(g)|^{p}\psi(x,g)\right)
^{1/p}=0\ .
$$
Using the triangle inequality we have
$$
|c-c_{n}|_{p,\psi,x}\leq |(a-a_{n})b|_{p,\psi,x}+|a_{n}(b-b_{n})|_{p,\psi,x}
:=I+II\ .
$$
According to (\ref{double}) for any $\epsilon>0$ we can decompose
$b$ in the sum $b'+b''$ where $b'=0$ outside a finite subset 
$S_{\epsilon}\subset X_{H}$ and $b''=0$ on $S_{\epsilon}$ such that
$|b''(g)|<\epsilon$ for all $g$. Note also that $|b'(g)|\leq 1$ for all $g$.
Also, (\ref{general}) and uniform continuity of $\log\psi$ 
on the compact set $bD$ imply that 
$|a-a_{n}|_{p,\psi,x}\leq km$ for some $k>0$. Finally, by continuity of
$f$ on $bD'$ we can find a number $N$ such that for any $n\geq N$ we
have $|(a-a_{n})\chi_{\epsilon}|_{p,\psi,x}<\epsilon$, where 
$\chi_{\epsilon}$ is the characteristic function of $S_{\epsilon}$.
Using all these facts we get for $n\geq N$
$$
I\leq |(a-a_{n})\chi_{\epsilon}b'|_{p,\psi,x}+|(a-a_{n})b''|_{p,\psi,x}
\leq\epsilon+km\epsilon=(1+km)\epsilon\ .
$$
To estimate $II$ observe that from (\ref{general}) and uniform continuity of 
$\log\psi$ follow that for  each $g\in X_{H}$
$$
(|a_{n}(g)|^{p}\psi(x,g))^{1/p}\leq |a_{n}|_{p,\psi,x}\leq 
k'|a_{n}|_{p,\psi,x_{n}}
\leq k' m\ 
$$
(for some $k'$). Moreover, since $F_{c_{1},c_{2}}|_{bD'}\in C_{p,1}(bD')$,
there is an integer $N'$ such that for any $n\geq N'$ we have
$|b-b_{n}|_{p,1,x}<\epsilon$. These two inequalities yield for
$n\geq N'$
$$
II\leq\sup_{g\in X_{H}}(|a_{n}(g)|^{p}\psi(x,g))^{1/p}|b-b_{n}|_{p,1,x}
\leq k'm\epsilon\ .
$$
Combining the estimates for $I$ and $II$ we obtain that 
$$
\lim_{n\to\infty}|c-c_{n}|_{p,\psi,x}=0\ .
$$
This is equivalent to $C_{p,\psi}$-continuity of $f_{c_{1},c_{2}}$ over
$x$.

Similarly one can check that if $p=\infty$ then $f_{c_{1},c_{2}}$ belongs to
$C_{1,1}(bD')$. We leave it as an excercise to the readers.\ \ \ \ \ $\Box$

Let us finish the proof of the theorem. According to Lemma \ref{smooth}
and the case (A) there is a function 
$f_{c_{1},c_{2}}'\in {\cal O}(D')\cap C(\overline{D'})$ such that
$f_{c_{1},c_{2}}'|_{bD'}=f_{c_{1},c_{2}}$. Note also that since
$|F_{c_{1},c_{2}}|\leq 1$ inequality (\ref{norm3}) yields
$$
|f_{c_{1},c_{2}}'|_{p,\psi}^{D}\leq c\sup_{z\in bD}|f|_{p,\psi,z}:=cm\ .
$$ 
Consider $f':=f_{c_{1},c_{2}}/F_{c_{1},c_{2}}$. Then 
$f'\in {\cal O}(D')\cap C(\overline{D'})$ and $f'|_{bD'}=f$. The uniqueness
property for holomorphic functions implies that $f'$ does not
depend on $c_{1}$ and $c_{2}$. Since $F_{c_{1},c_{2}}$ converges uniformly 
on compact subsets of $M'$ to $1$ as $c_{2}\to 0$ from the last inequality 
we get
$$
|f'|_{p,\psi}^{D}\leq c\sup_{z\in bD}|f|_{p,\psi,z}\ .
$$
Therefore $f'\in {\cal H}_{p,\psi}(D')\cap C(\overline{D'})$.

The proof of the theorem for $s=0$ is complete. If, in addition,
$f\in C^{s}(bD)$ for $1\leq s\leq k$, then in fact the extended function
$f'\in C^{s}(\overline{D'})$ (see, e.g., Theorem 3.14 in [Ha], and
the discussion that follows it).\ \ \ \ \ $\Box$\\
{\bf Proof of Corollary \ref{co2}.}
Let us consider the function $F_{c_{1},c_{2}}$ from (\ref{double}).
Suppose that $f\in {\cal O}(U')$ satisfies the hypotheses of Corollary 
\ref{co2} with $c_{1}$, $c_{2}$ and $c$ as in the definition of
$F_{c_{1},c_{2}}$. Then $f_{c_{1},c_{2}}:=fF_{c_{1},c_{2}}\in 
H^{\infty}(U')$ with the norm bounded by 1. Further, by the hypotheses we can 
find a connected $C^{\infty}$ compact submanifold $bS\subset U$ that bounds 
a domain $S$ containing $K$.
By Theorem \ref{te4} (applied to $S':=r^{-1}(S)$ and $bS':=r^{-1}(bS)$) the 
function $f_{c_{1},c_{2}}$ admits an extension 
$f'_{c_{1},c_{2}}\in H^{\infty}(D')$. Since 
$f'_{c_{1},c_{2}}=f_{c_{1},c_{2}}$ on $U'$ and
$h(z):=\sup_{y\in r^{-1}(z)}|f'_{c_{1},c_{2}}(y)|$ is a continuous 
plurisubharmonic function on $D$,
$$
|f'_{c_{1},c_{2}}|_{\infty}^{D}=\sup_{z\in U'}|f(z)|\leq 1\ .
$$
Then the function $f':=f'_{c_{1},c_{2}}/F_{c_{1},c_{2}}\in {\cal O}(D')$
extends $f$ and  satisfies
$$
|f'(z)|\leq e^{c_{3}e^{c_{1}d_{o}(z)}}\ ,\ \ \ z\in D'\ ,
$$
with $c_{3}$ as in (\ref{double}).\ \ \ \ \ $\Box$ 
\sect{\hspace*{-1em}. Proof of Theorem \ref{te5}.}
Let $r:M_{G}\to M$ be the regular covering of $M$ satisfying condition
(\ref{e1}) (for some $\widetilde M$ and $N$) with transformation group $G$. 
Let $G_{1}\subset G$ be a 
subgroup of a finite index and let $r_{1}:M_{1}\to M$ be the covering 
with fibre $G/G_{1}$. Then there are coverings $N_{1}$ and 
$\widetilde M_{1}$ of $N$ and $\widetilde M$ with 
fibre $G/G_{1}$ such that
$M_{1}\subset\subset\widetilde M_{1}\subset N_{1}$. Clearly this triple
satisfies condition (\ref{e1}), as well. Thus without loss of generality
we may assume that $M:=M_{1}$, $G:=G_{1}$, $\widetilde M:=\widetilde M_{1}$
and $N:=N_{1}$, and so $G\in {\cal T}(M)$. The latter means that $G$ admits
a linear representation $\rho$ into $GL_{k}(\Co)$ and that the flat vector
bundle $E_{\rho}$ on $N$ associated with $\rho$ is topologically trivial. 
Then the restriction $E_{\rho}|_{\widetilde M}$ is
topologically trivial. Since $\widetilde M$ is Stein,
according to the Oka-Grauert principle (see [G2]),
$E_{\rho}|_{\widetilde M}$ is holomorphically trivial. In particular,
$\rho$ can be obtained as the monodromy of an equation
$dF=\omega F$ on $\widetilde M$ where $\omega$ is a matrix-valued holomorphic
1-form on $\widetilde M$ satisfying $d\omega-\omega\wedge\omega=0$.
Let $\widetilde\omega=r^{*}\omega$ be the pullback of $\omega$ on
$\widetilde M_{G}=r^{-1}(\widetilde M)$. (Here $r:N_{G}\to N$ is the regular
covering of $N$ with the transformation group $G$ so that
$M_{G}$ and $\widetilde M_{G}$ are domains in $N_{G}$.) 
Then there exists a function
$\widetilde F\in {\cal O}(\widetilde M_{G},GL_{k}(\Co))$ such that
$d\widetilde F=\widetilde\omega\widetilde F$. This follows from the fact
that the monodromy of the last equation is the restriction of $\rho$
to $\pi_{1}(M_{G})$ and so it is trivial (since 
$\pi_{1}(M_{G})\subset Ker\ \!\rho$).
Note that $\widetilde F$ can be obtained 
by Picard iteration applied to $\widetilde\omega$. Since
$\overline{M}$ is a compact subset of $\widetilde M$ (and so
$\omega|_{M}$ is bounded), the Picard iteration produces for some 
positive $c=c(M,\omega)$ the estimate 
\begin{equation}\label{gro}
||\widetilde F(z)||_{2}\leq e^{cd_{o}(z)}\ ,\ \ \ z\in M_{G}\ ,
\end{equation}
where $||\cdot||_{2}$ is the 
$l_{2}$-norm on $GL_{k}(\Co)$ and $o\in M_{G}$. (Here as before
$d_{o}$ is the distance from $o$ in the path metric induced by a 
Riemannian metric pulled back from $N$.) Moreover, for every
$z\in\widetilde M_{G}$ there exists a matrix $C_{z}\in GL_{k}(\Co)$ such
that $\widetilde F(gz)=C_{z}^{-1}\rho(g)C_{z}$ for any $g\in G$.
These are standard facts of the theory of flat connections. In particular,
from the last identity we derive easily that
$\widetilde F$ separates points in every orbit of the action of $G$ on 
$\overline{M}_{G}$. 

Next, let $f$ be the function from Corollary \ref{co1}. 
Then by (\ref{gro}) we get (for some $c_{1}=c_{1}(f,\alpha)$)
$$
||e^{-\alpha f(z)}\widetilde F(z)||_{2}\leq c_{1}e^{(c-\alpha)d_{o}(z)}\ ,
\ \ \ z\in M_{G}\ .
$$
From here arguing as in the proof of Corollary \ref{co1} we deduce that
for a sufficiently big $\alpha$ all entries of the matrix
$e^{-\alpha f}\cdot\widetilde F$ belong to 
${\cal H}_{2,1}(M_{G})\cap C(\overline{M}_{G})$.
Now the family consisting of these entries and the function
$e^{-\alpha f}$ separate all points in any orbit of the action of $G$ on
$\overline{M}_{G}$. For otherwise, there are $x,y\in\overline{M}_{G}$, $y=gx$,
$g\neq 1$, $g\in G$, such that 
$e^{-\alpha f(x)}\widetilde F(x)=e^{-\alpha f(y)}\widetilde F(y)$ and
$e^{-\alpha f(x)}=e^{-\alpha f(y)}$. But this implies that 
$\widetilde F(x)=\widetilde F(y)$, a contradiction. Finally,
since $M\subset\subset\widetilde M$ and $\widetilde M$ is Stein, by the
Remmert embedding theorem there are holomorphic functions 
$h_{1},\dots, h_{l}$ from $H^{\infty}(M)\cap C(\overline{M})$ that separate
all points in $\overline{M}$. We set
$\widetilde h_{i}=e^{-\alpha f}r^{*}h_{i}$, $1\leq i\leq l$. Then by the
definition $\widetilde h_{i}\in {\cal H}_{2,1}(M_{G})\cap C(\overline{M}_{G})$
and so the family consisting of all $\widetilde h_{i}$, 
entries of $e^{-\alpha f}\widetilde F$ and $e^{-\alpha f}$
separates all points in $\overline{M}_{G}$.

The proof of Theorem \ref{te5} is complete.\ \ \ \ \ $\Box$


\end{document}